\documentclass[a4paper]{article}
\usepackage[abbrvbib,preprint]{jmlr2e}

\usepackage[english]{babel}
\usepackage[utf8x]{inputenc}
\usepackage[T1]{fontenc}

\usepackage{amsmath,amssymb}
\usepackage{graphicx}
\usepackage[font=small,labelfont=bf]{caption}

\usepackage[shortlabels]{enumitem}

\ShortHeadings{Implicit ridge regularization}{Kobak, Lomond, Sanchez}

\title{Optimal ridge penalty for real-world high-dimensional data can be zero or negative due to the implicit ridge regularization}

\author{\name Dmitry Kobak \email dmitry.kobak@uni-tuebingen.de \\
       \addr Institute for Ophthalmic Research \\University of T\"ubingen\\
       Otfried-M\"uller-Stra{\ss}e 25, 72076 T\"ubingen
       \AND
       \name Jonathan Lomond \\
       \addr Toronto, Canada
       \AND Benoit Sanchez\\
       \addr SmartAdServer, Paris, France}

\DeclareMathOperator*{\argmin}{arg\,min}

\newcommand{\X}{\mathbf X}
\newcommand{\y}{\mathbf y}
\newcommand{\x}{\mathbf x}
\newcommand{\I}{\mathbf I}
\newcommand{\bbeta}{{\boldsymbol\beta}}
\newcommand{\bbetahat}{\boldsymbol{\hat\beta}}
\renewcommand{\S}{\boldsymbol\Sigma}

\begin{document}
\maketitle

\begin{abstract}
A conventional wisdom in statistical learning is that large models require strong regularization to prevent overfitting. 
Here we show that this rule can be violated by  linear regression in the underdetermined $n\ll p$ situation under realistic conditions.
Using simulations and real-life high-dimensional data sets, we demonstrate that an explicit positive ridge penalty can fail to provide any improvement over the minimum-norm least squares estimator. Moreover, the optimal value of ridge penalty in this situation can be negative. This happens when the high-variance directions in the predictor space can predict the response variable, which is often the case in the real-world high-dimensional data. In this regime, low-variance directions provide an implicit ridge regularization and can make any further positive ridge penalty detrimental. We prove that augmenting any linear model with random covariates and using minimum-norm estimator is asymptotically equivalent to adding the ridge penalty. We use a spiked covariance model as an analytically tractable example and prove that the optimal ridge penalty in this case is negative when $n\ll p$.
\end{abstract}


\section{Introduction}

In recent years, there has been increasing interest in prediction problems in which the sample size $n$ is much smaller than the dimensionality of the data $p$. This situation is known as $n\ll p$ and often arises in computational chemistry and biology, e.g. in chemometrics, brain imaging, or genomics \citep{hastie2009elements}. The standard approach to such problems is ``to bet on sparsity'' \citep{hastie2015statistical} and to use linear models with regularization performing feature selection, such as the lasso \citep{tibshirani1996regression}, the elastic net \citep{zou2005regularization}, or the Dantzig selector \citep{candes2007dantzig}.

In this paper we discuss ordinary least squares (OLS) linear regression with loss function \begin{equation}\mathcal L = \lVert \y - \X \bbeta\rVert^2,\end{equation} where $\X$ is a $n\times p$ matrix of predictors and $\y$ is a $n\times 1$ matrix of responses. Assuming $n>p$ and full-rank $\X$, the unique solution minimizing this loss function is given by  \begin{equation}\bbetahat_\mathrm{OLS} = (\X^\top \X)^{-1} \X^\top \y.\end{equation} This estimator is unbiased and has small variance when $n\gg p$. As $p$ grows for a fixed $n$, $\X^\top \X$ becomes poorly conditioned, increasing the variance and leading to overfitting. The expected error can be decreased by shrinkage as provided e.g. by the ridge estimator \citep{hoerl1970ridge}, a special case of Tikhonov regularization \citep{tikhonov1963solution}, \begin{equation}\bbetahat_\lambda = (\X^\top \X + \lambda \I)^{-1} \X^\top \y,\end{equation} which minimizes the loss function with an added $\ell_2$ penalty \begin{equation}\mathcal L_\lambda = \lVert \y - \X \bbeta\rVert^2 + \lambda\lVert\bbeta\rVert^2.\end{equation}

The closer $p$ is to $n$, the stronger the overfitting and the more important it is to use regularization. It seems intuitive that when $p$ becomes larger than $n$, regularization becomes indispensable and small values of $\lambda \approx 0$ would yield hopeless overfitting. A popular textbook \citep{james2013introduction}, for example, claims that \textit{``though it is possible to perfectly fit the training data in the high-dimensional setting, the resulting linear model will perform extremely poorly on an independent test set, and therefore does not constitute a useful model.''} Here we show that this intuition is incomplete.

Specifically, we empirically demonstrate and mathematically prove the following: 
\begin{enumerate}[(i)]
    \item when $n\ll p$, the $\lambda\to 0$ limit, corresponding to the minimum-norm OLS solution, can have good generalization performance; 
    \item additional ridge regularization with $\lambda>0$ can fail to provide any further improvement; 
    \item moreover, the optimal value of $\lambda$ in this regime can be \textit{negative}; 
    \item this happens when response variable is predicted by the high-variance directions while the low-variance directions together with the minimum-norm requirement effectively perform shrinkage and provide implicit ridge regularization.
\end{enumerate}

Our results provide a simple counter-example to the common understanding that large models with little regularization do not generalize well. This has been  pointed out as a puzzling property of deep neural networks \citep{zhang2016understanding}, and has been subject to a very active ongoing research since then, performed independently from our work (the first version of this manuscript was released as a preprint in May 2018). Several groups reported that very different statistical models can display {\textbackslash}/{\textbackslash}-shaped (\textit{double descent}) risk curves as a function of model complexity, extending the classical U-shaped risk curves and having small or even the smallest risk in the $p\gg n$ regime \citep{advani2017high, belkin2019reconciling, spigler2019jamming}. The same phenomenon was later demonstrated for modern deep learning architectures \citep{nakkiran2020deep}. In the context of linear or kernel methods, the high-dimensional regime when the model is rich enough to fit any training data with zero loss, has been called \textit{ridgeless} regression or  \textit{interpolation} \citep{liang2018just, hastie2019surprises}. The fact that such interpolating estimators can have low risk has been called \textit{benign overfitting} \citep{bartlett2019benign, chinot2020benign} or \textit{harmless interpolation} \citep{muthukumar2019harmless}.

Our finding (i) is in line with this body of parallel literature. Findings (ii) and (iii) have not, to the best of our knowledge, been described anywhere else. Existing studies of high-dimensional ridge regression found that, under some generic assumptions, the ridge risk at some $\lambda>0$ always dominates the minimum-norm OLS risk \citep{dobriban2018high, hastie2019surprises}. Our results highlight that the optimal value of ridge penalty can be zero or even negative, suggesting that real-world $n\ll p$ data sets can have very different statistical structure compared to the common theoretical models \citep{dobriban2018high}. Finding (ii) has been observed for kernel methods \citep{liang2018just} and for random features regression \citep{mei2019generalization}; our results demonstrate that (ii) can happen in a simpler situation of ridge regression with Gaussian features. We are not aware of any existing work reporting that the optimal ridge penalty can be \textit{negative}, as per our finding (iii). Finally, finding (iv) is related to the results of \citet{bibas2019new} and \citet{bartlett2019benign}; the connection between the minimum-norm OLS and the ridge estimators was also studied by \citet{derezinski2019exact}.


The code in Python can be found at \url{http://github.com/dkobak/high-dim-ridge}.

\section{Results}

\subsection{A case study of ridge regression in high dimensions}

We used the \texttt{liver.toxicity} dataset \citep{bushel2007simultaneous} from the R package \texttt{mixOmics} \citep{rohart2017mixomics} as a motivational example to demonstrate the phenomenon. This dataset contains microarray expression levels of $p=3116$ genes and 10 clinical chemistry measurements in liver tissue of $n=64$ rats. We centered and standardized all the variables before the analysis.

\begin{figure}
\includegraphics[width=\textwidth]{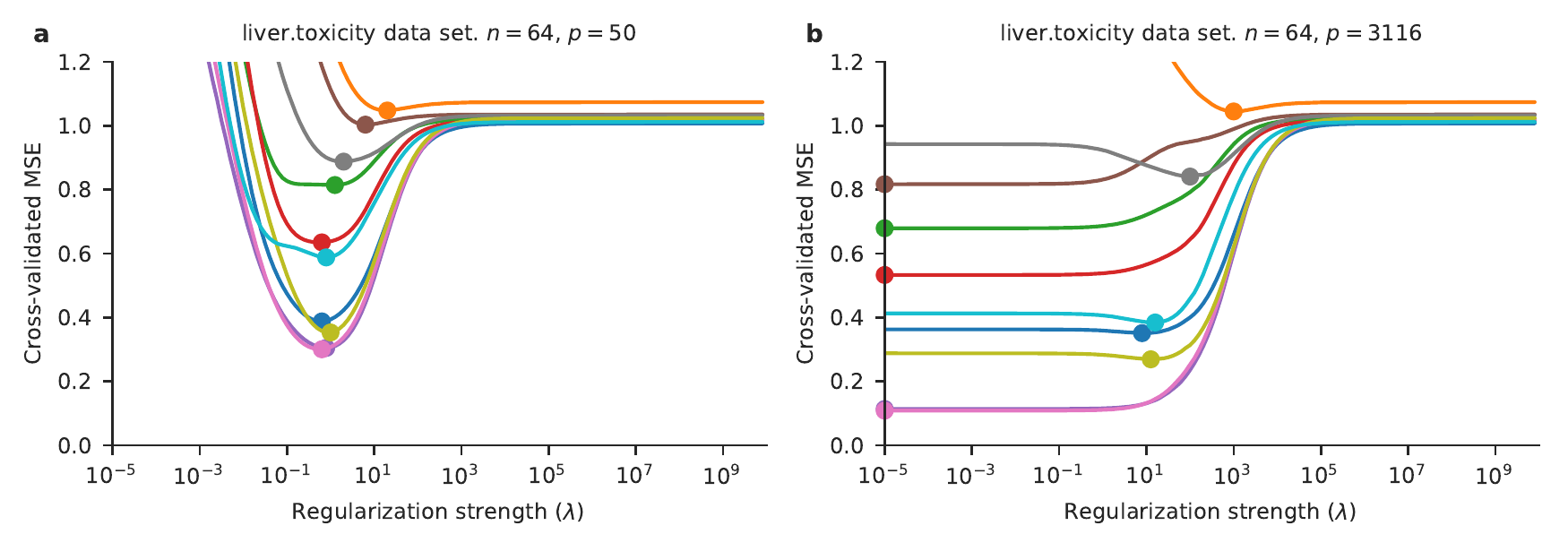}
\caption{Cross-validation estimate of ridge regression performance for the \texttt{liver.toxicity} dataset. \textbf{a.}~Using $p=50$ random predictors. \textbf{b.}~Using all $p=3116$ predictors. Lines correspond to 10 dependent variables. Dots show minimum values.}
\label{fig:casestudy}
\end{figure}

We used \texttt{glmnet} library \citep{friedman2010regularization} to predict each chemical measurement from the gene expression data using ridge regression. \texttt{Glmnet} performed 10-fold cross-validation (CV) for various values of regularization parameter $\lambda$. 
We ran CV separately for each of the 10 dependent variables. When we used $p=50$ random predictors, there was a clear minimum of mean squared error (MSE) for some $\lambda_\mathrm{opt}>0$, and smaller values of $\lambda$ yielded much higher MSE, i.e. led to overfitting (Figure~\ref{fig:casestudy}a). This is in agreement with \citet{hoerl1970ridge} who proved  that when $n<p$, the optimal penalty $\lambda_\mathrm{opt}$ is always larger than zero. The CV curves had a similar shape when $p\gtrsim n$, e.g. $p=75$.

However, when we used all $p\gg n$ predictors, the curves changed dramatically (Figure~\ref{fig:casestudy}b). For five dependent variables out of ten, the lowest MSE corresponded to the smallest value of $\lambda$ that we tried. Four other dependent variables had a minimum in the middle of the $\lambda$ range, but the limiting MSE value at $\lambda\to 0$ was close to the minimal one. This is counter-intuitive: despite having more predictors than samples, tiny values of $\lambda \approx 0$ provide optimal or near-optimal estimator.

We observed the same effect in various other genomics datasets with $n\ll p$ \citep{kobak2018sparse}. We believe it is a general phenomenon and not a peculiarity of this particular dataset.

\subsection{Minimum-norm OLS estimator}

When $n<p$, the limiting value of the ridge estimator at $\lambda \to 0$ is the minimum-norm OLS estimator. It can be shown using a thin singular value decomposition (SVD) of the predictor matrix $\X = \mathbf{USV}^\top$ (with $\mathbf S$ square and all its diagonal values non-zero):
\begin{equation} \label{eq:minorm}
\bbetahat_0 = \lim_{\lambda\to 0}\bbetahat_\lambda = \lim_{\lambda\to 0}(\X^\top \X + \lambda \I)^{-1} \X^\top \y = \lim_{\lambda\to 0}\mathbf V \frac{\mathbf S}{\mathbf S^2+\lambda}\mathbf U^\top \y = \mathbf V \mathbf S^{-1}\mathbf U^\top \y = \X^+ \y, 
\end{equation}
where $\X^+=\X^\top(\X\X^\top)^{-1}$ denotes pseudo-inverse of $\X$ and operations on the diagonal matrix $\mathbf S$ are assumed to be element-wise and applied only to the diagonal. 

The estimator $\bbetahat_0$ gives one possible solution to the OLS problem and, as any other solution, it provides a perfect fit on the training set:
\begin{equation}\lVert \y - \X \bbetahat_0\rVert^2 = \lVert \y - \X \X^+ \y\rVert^2 = \lVert \y - \y\rVert^2 = 0.\end{equation}
The $\bbetahat_0$ solution is the one with minimum $\ell_2$ norm:
\begin{equation}\bbetahat_0 = \argmin\Big\{\lVert\bbeta\rVert^2 \Bigm| \lVert\y - \X \bbeta\rVert^2 = 0\Big\}.\end{equation}
Indeed, any other solution can be written as a sum of $\bbetahat_0$ and a vector from the $(p-n)$-dimensional subspace orthogonal to the column space of $\mathbf V$. Any such vector yields a valid OLS solution but increases its norm compared to $\bbetahat_0$ alone.

This allows us to rephrase the observations made in the previous section as follows: when $n\ll p$, the minimum-norm OLS estimator can be better than any ridge estimator with $\lambda>0$.

\subsection{Simulation using spiked covariance model}
\label{section:model}

We qualitatively replicated this empirically observed phenomenon with a simple model where all $p$ predictors are positively correlated to each other and all have the same effect on the response variable.

Let $\mathbf x\sim\mathcal N(0,\boldsymbol\Sigma)$ be a $p$-dimensional vector of predictors with covariance matrix $\boldsymbol\Sigma$ having all diagonal values equal to $1+\rho$ and all non-diagonal values equal to $\rho$. This is known as \textit{spiked covariance model}: $\boldsymbol\Sigma = \mathbf I + \rho \mathbf 1 \mathbf 1^\top$ deviates from the spherical covariance $\mathbf I$ only in one dimension. Let the response variable be $y=\mathbf x^\top \boldsymbol \beta + \varepsilon$, where $\varepsilon\sim\mathcal N(0,\sigma^2)$ and $\boldsymbol \beta=(b,b\ldots b)^\top$ has all identical elements. We select $b=\sigma \sqrt{\alpha/(p+p^2\rho)}$ in order to achieve signal-to-noise ratio $\operatorname{Var}[\mathbf x^\top \boldsymbol \beta]/\operatorname{Var}[\varepsilon]=\operatorname{Var}[\mathbf x^\top \boldsymbol \beta]=\alpha$. In all simulations we fix $\sigma^2=1$, $\rho=0.1$ and $\alpha=10$.

\begin{figure}[t]
\includegraphics[width=\textwidth]{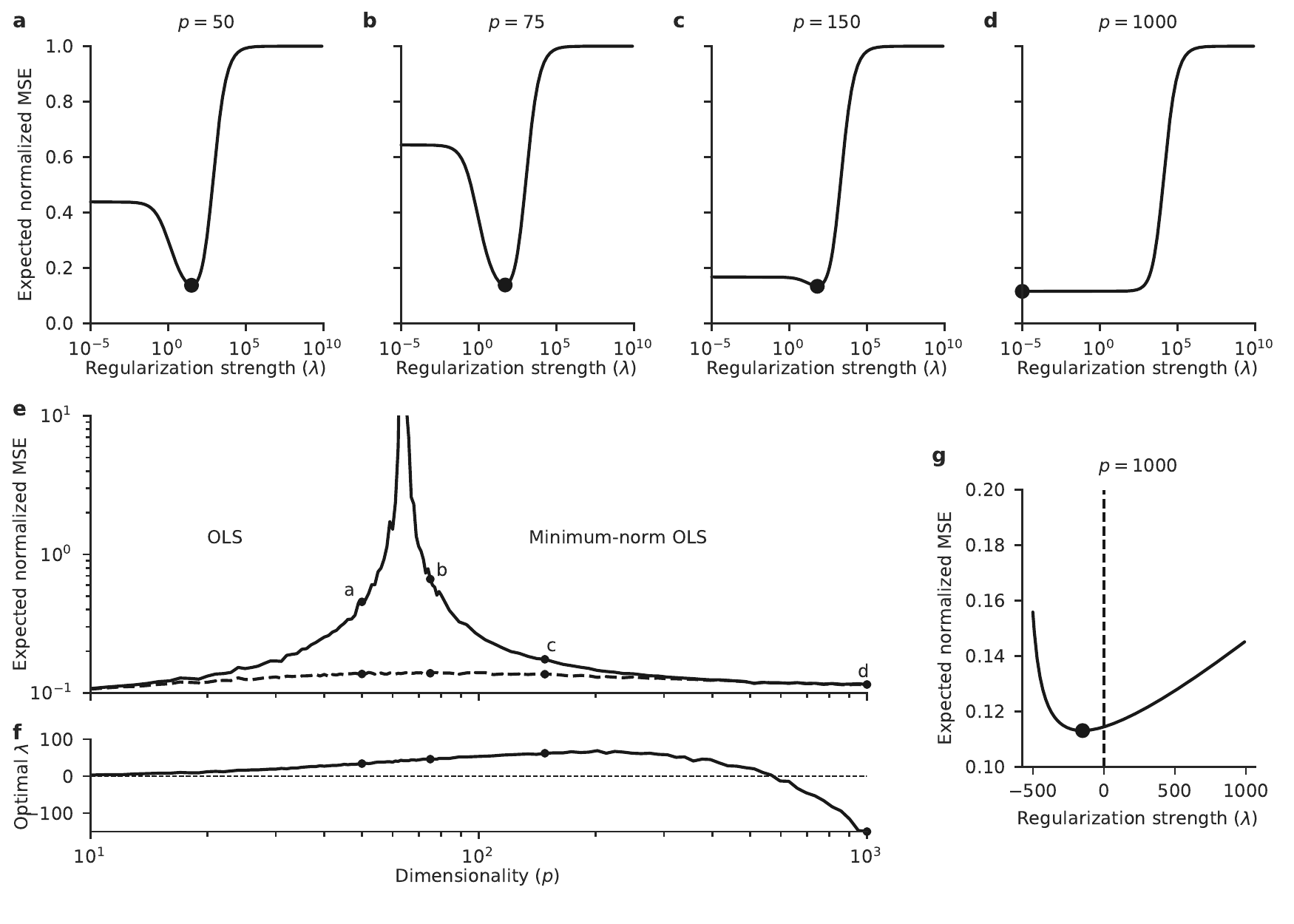}
\caption{\textbf{a--d.} Expected normalized MSE of ridge estimators using a model with correlated predictors. On all subplots $n=64$. Subplots correspond to the number of predictors $p$ taking values 50, 75, 150, and 1000. Dots mark the points with minimum risk.
\textbf{e.} Expected normalized MSE of OLS (for $n<p$) and minimum-norm OLS (for $p>n$) estimators using the same model with $p\in[10,1000]$. Dots mark the dimensionalities corresponding to subplots (a--d). Dashed line: the expected normalized MSE of the optimal ridge estimator. \textbf{f.} The values of $\lambda$ minimizing the expected risk. For $p\gtrsim 600$, the optimal value of ridge penalty was negative: $\lambda_\mathrm{opt}<0$. \textbf{f.} Expected normalized MSE of ridge estimators for $p=1000$ including negative values of $\lambda$. The minimum was attained at $\lambda_\mathrm{opt}=-150$.}
\label{fig:model}
\end{figure}

Using this model with different values of $p$, we generated many ($N_\mathrm{rep}=100$) training sets $(\mathbf X, \mathbf y)$ with $n=64$ each, as in the \texttt{liver.toxicity} dataset analyzed above. Using each training set, we computed $\bbetahat_\lambda=\mathbf V \frac{\mathbf S}{\mathbf S^2+\lambda}\mathbf U^\top \y$ for various values of $\lambda$ and then found MSE (risk) of $\bbetahat_\lambda$ using the formula
\begin{equation}
    R(\bbetahat_\lambda) = \mathbb E_{\mathbf x, \varepsilon}\big[\big((\mathbf x^\top\bbeta + \varepsilon) - \mathbf x^\top\bbetahat_\lambda\big)^2\big] = (\bbetahat_\lambda-\bbeta)^\top\boldsymbol\Sigma(\bbetahat_\lambda-\bbeta) + \sigma^2.
\end{equation} We normalized the MSE by $\operatorname{Var}[y]=\bbeta^\top\boldsymbol\Sigma\bbeta+\sigma^2=(\alpha+1)\sigma^2$. Then we averaged normalized MSEs across $N_\mathrm{rep}$ training sets to get an estimate of the expected normalized MSE. The results for $p\in\{50,75,150,1000\}$ (Figure~\ref{fig:model}a--d) match well to what we previously observed in real data (Figure~\ref{fig:casestudy}): when $n>p$ or $n\lesssim p$, the MSE had a clear minimum for some positive value of $\lambda$. But when $n\ll p$, the minimum MSE was achieved by the $\lambda=0$ minimum-norm OLS estimator.

Figure~\ref{fig:model}e shows the expected normalized MSE of the OLS and the minimum-norm OLS estimators for $p\in[10,1000]$. The true signal-to-noise ratio was always $\alpha=10$, so the best attainable normalized MSE was always $1/(10+1)\approx 0.09$. With $p=10$, OLS yielded a near-optimal performance. As $p$ increased, OLS began to overfit and each additional predictor increased the MSE. Near $p\approx n=64$ the expected MSE became very large, but as $p$ increased even further, the MSE of the minimum-norm OLS quickly decreased again.

The risk of the optimal ridge estimator was close to the oracle risk for all dimensionalities (Figure~\ref{fig:model}e, dashed line), and did not show any divergence at $p=n$. However, as $p>n$ grew, the gain compared to the minimum-norm OLS estimator became smaller and smaller and in the $p\gg n$ regime eventually disappeared. Moreover, for sufficiently large values of $p$, the optimal regularization value $\lambda_\mathrm{opt}$ became negative (Figure~\ref{fig:model}f). We found it to be the case for $p\gtrsim 600$. In sufficiently large dimensionalities, the expected risk as a function of $\lambda$ had a minimum not at zero (Figure~\ref{fig:model}d), but at some negative value of $\lambda$ (Figure~\ref{fig:model}g). For $p=1000$, the lowest risk was achieved at $\lambda_\mathrm{opt}=-150$ (Figure~\ref{fig:model}g).

\begin{figure}[t]
\includegraphics[width=\textwidth]{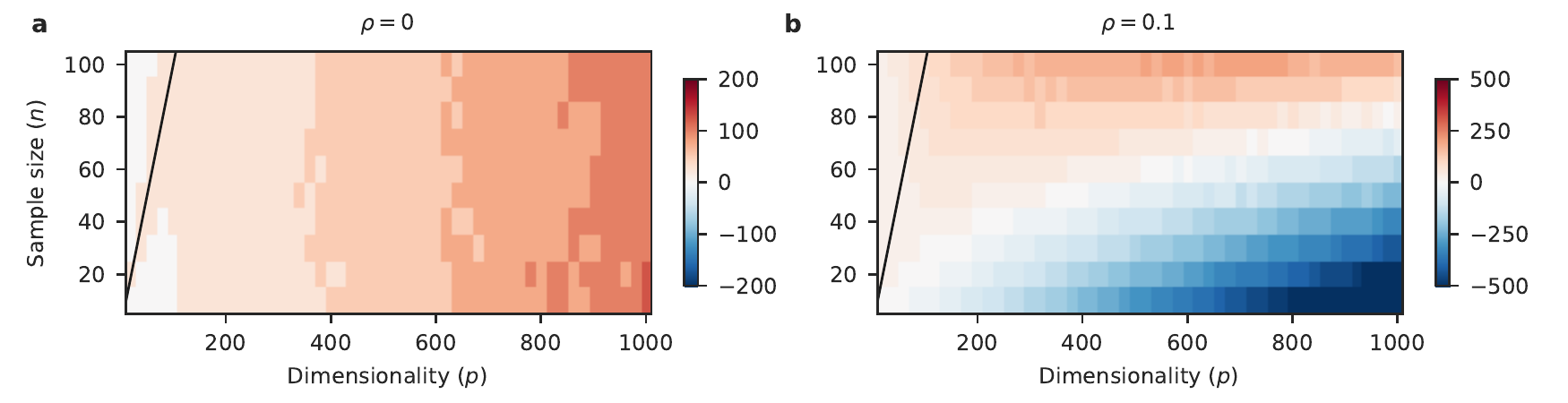}
\caption{\textbf{a.} The optimal regularization parameter $\lambda_\mathrm{opt}$ as a function of sample size ($n$) and dimensionality ($p$) in the model with uncorrelated predictors ($\rho=0$). In this case $\lambda_\mathrm{opt}=p\sigma^2/\lVert\boldsymbol\beta\rVert=p/\alpha$. Black line corresponds to $n=p$.
\textbf{b.} The optimal regularization parameter $\lambda_\mathrm{opt}$ in the model with correlated predictors ($\rho=0.1$).}
\label{fig:optlambdas}
\end{figure}

To investigate this further, we found the optimal regularization value $\lambda_\mathrm{opt}$ for different sample sizes $n\in[10,100]$ and different dimensionalities $p\in[20,1000]$ (Figure~\ref{fig:optlambdas}). For the spherical covariance matrix ($\rho=0$), $\lambda_\mathrm{opt}$ did not depend on the sample sizes and grew linearly with dimensionality (Figure~\ref{fig:optlambdas}a), in agreement with the analytical formula $\lambda_\mathrm{opt}=p\sigma^2/\lVert\boldsymbol\beta\rVert^2=p/\alpha$ \citep{nakkiran2020optimal}. But in our model with $\rho=0.1$, the optimal value $\lambda_\mathrm{opt}$ in sufficiently high dimensionality was negative for any sample size. The smallest dimensionality necessary for this to happen grew with the sample size (Figure~\ref{fig:optlambdas}b).

This result might appear to contradict the literature; for example, \citet{dobriban2018high} and later \citet{hastie2019surprises} studied high-dimensional asymptotics of ridge regression performance for $p,n\to\infty$ while $p/n=\gamma$ and proved, among other things, that the optimal $\lambda$ is always positive. Their results hold for an arbitrary covariance matrix $\boldsymbol\Sigma$ when the elements of $\boldsymbol\beta$ are random with mean zero. The key property of our simulation is that $\boldsymbol\beta$ is not random and does not point in a random direction; instead, it is aligned with the first principal component (PC1) of $\boldsymbol\Sigma$.

While such a perfect alignment can never hold exactly in real-world data, it is plausible that $\boldsymbol\beta$ often points in a direction of sufficiently high predictor variance. Indeed, principal component regression (PCR) that discards all low-variance PCs and only uses high-variance PCs for prediction is known to work well for many real-world $n\ll p$ data sets \citep{hastie2009elements}. In the next section we show that the low-variance PCs can provide an implicit ridge regularization.

\subsection{Implicit ridge regularization provided by random low-variance predictors}

Here we prove that augmenting a model with randomly generated low-variance predictors is  is asymptotically equivalent to the ridge shrinkage.

\begin{theorem} Let $\bbetahat_\lambda$ be a ridge estimator of $\bbeta\in\mathbb R^p$ in a linear model $y= \x^\top \bbeta + \varepsilon$, given some training data $(\X, \y)$ and some value of $\lambda$. We construct a new estimator $\bbetahat_q$ by augmenting $\X$ with $q$ columns $\X_q$ with i.i.d. elements, randomly generated with mean $0$ and variance $\lambda/q$, fitting the model with minimum-norm OLS, and taking only the first $p$ elements. Then \begin{equation*}\bbetahat_q \xrightarrow[q\to\infty]{\mathrm{a.s.}} \bbetahat_\lambda.\end{equation*} 
In addition, for any given $\x$, let $\hat y_\lambda=\x^\top \bbetahat_\lambda$ be the response predicted by the ridge estimator, and $\hat y_\mathrm{augm}$ be the response predicted by the augmented model including the additional $q$ parameters using $\x$ extended with $q$ random elements (as above). Then:
\begin{equation*}\hat y_\mathrm{augm} \xrightarrow[q\to\infty]{\mathrm{a.s.}} \hat y_\lambda.\end{equation*}
\end{theorem}

\newcommand{\Xaugm}{\mathbf X_\mathrm{augm}}
\newcommand{\betaaugm}{\boldsymbol{\hat\beta}_\mathrm{augm}}
\newcommand{\Xq}{\mathbf X_q}

\begin{proof}
Let us write $\Xaugm = \begin{bmatrix}\X & \Xq\end{bmatrix}$. The minimum-norm OLS estimator can be written as \begin{equation}\betaaugm = \Xaugm^+\y = \Xaugm^\top (\Xaugm \Xaugm^\top)^{-1} \y.\end{equation}
By the strong law of large numbers,
\begin{equation}\Xaugm \Xaugm^\top = \X\X^\top + \Xq \Xq^\top \to \X\X^\top + \lambda \I_n.\end{equation}
The first $p$ components of $\betaaugm$ are
\begin{equation}\bbetahat_q=\X^\top(\Xaugm \Xaugm^\top)^{-1}\y\to \X^\top ( \X\X^\top + \lambda \I_n)^{-1} \y.\end{equation} 
Note that
$(\X^\top \X + \lambda \I_p)\X^\top=\X^\top(\X\X^\top + \lambda \I_n)$.
Multiplying this equality by $(\X^\top \X + \lambda \I_p)^{-1}$ on the left and $(\X\X^\top + \lambda \I_n)^{-1}$ on the right, we obtain the following standard identity:
\begin{equation}\X^\top(\X\X^\top + \lambda \I_n)^{-1}=(\X^\top \X + \lambda \I_p)^{-1}\X^\top.\end{equation}
Finally:
\begin{equation}\bbetahat_q\to(\X^\top \X + \lambda \I_p)^{-1}\X^\top \y = \bbetahat_\lambda.\end{equation} 

To prove the second statement of the Theorem, let us write $\x_\mathrm{augm}=\begin{bmatrix}\x\\ \x_q\end{bmatrix}$. The predicted value using the augmented model is:
\begin{align}
\hat y_\mathrm{augm} = \x_\mathrm{augm}^\top \bbetahat_\mathrm{augm} &= \x_\mathrm{augm}^\top \Xaugm^\top (\Xaugm \Xaugm^\top)^{-1} \y \\ 
&= \begin{bmatrix}\x\\\x_q\end{bmatrix}^\top \begin{bmatrix}\X & \X_q\end{bmatrix}^\top (\X\X^\top + \Xq \Xq^\top)^{-1} \y \\
&= \x^\top \X^\top (\X\X^\top + \Xq \Xq^\top)^{-1} \y + \x_q^\top \X_q^\top (\X\X^\top + \Xq \Xq^\top)^{-1} \y \\
&\to \x^\top \bbetahat_\lambda + 0_{1\times n} (\X\X^\top + \lambda \I_n)^{-1} \y \\
& = \x^\top \bbetahat_\lambda = \hat y_\lambda.
\end{align}
\end{proof}

Note that the Theorem requires the random predictors to be independent from each other, but does \textit{not} require them to be independent from the existing predictors or from the response variable. 

From the first statement of the Theorem it follows that the expected MSE of the truncated estimator $\bbetahat_q$ converges to the expected MSE of the ridge estimator $\bbetahat_\lambda$. From the second statement it follows that the expected MSE of the augmented estimator on the augmented data also converges to the expected MSE of the ridge estimator.

\begin{figure}
\includegraphics[width=\textwidth]{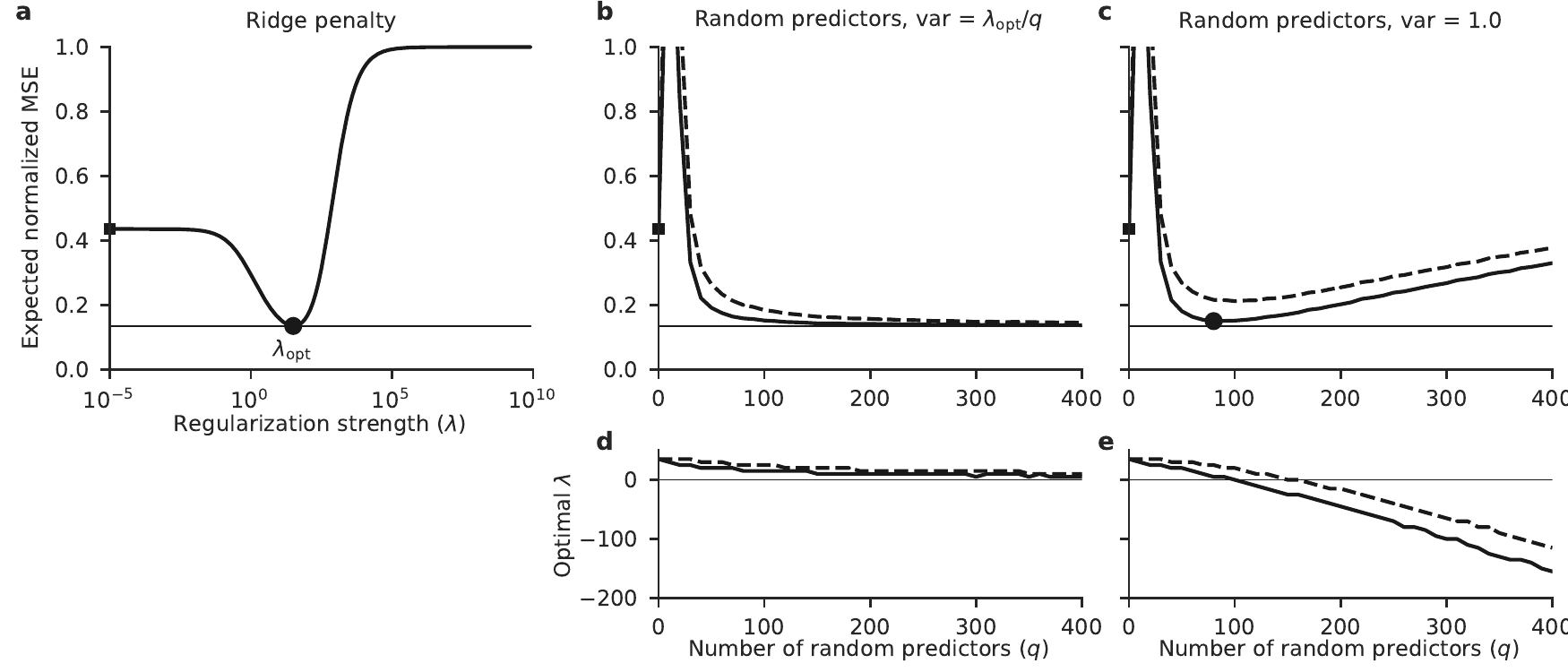}
\caption{\textbf{a.} Expected MSE as a function of ridge penalty in the toy model with $p=50$ weakly correlated predictors that are all weakly correlated with the response ($n=64$). This is the same plot as in Figure~\ref{fig:model}a. The dot denotes minimal risk and the square denotes the MSE of the OLS estimator ($\lambda=0$). The horizontal line shows the optimal risk corresponding to $\lambda_\mathrm{opt}$. \textbf{b.} Augmenting the model with up to $q=400$ random predictors with variance $\lambda_\mathrm{opt}/q$. Solid line corresponds to $\bbetahat_q$ (i.e. $\bbetahat_\mathrm{augm}$ truncated to $p$ predictors); dashed line corresponds to the full $\bbetahat_\mathrm{augm}$. \textbf{c.} Augmenting the model with up to $q=400$ random predictors with variance equal to 1. \textbf{d.} The optimal ridge penalty $\lambda_\mathrm{opt}$ in the model augmented with random predictors with adaptive variance, as in panel (b). \textbf{e.} The optimal ridge penalty $\lambda_\mathrm{opt}$ in the model augemented with random predictors with variance 1, as in panel (c).}
\label{fig:randompredictors}
\end{figure}

We extended the simulation from Section~\ref{section:model} to confirm this experimentally. We considered the same toy model as above with $n=64$ and $p=50$. Figure~\ref{fig:randompredictors}a (identical to Figure~\ref{fig:model}a) shows the expected MSE of ridge estimators for different values of $\lambda$. The optimal $\lambda$ in this case happened to be $\lambda_\mathrm{opt} = 31$. Figure~\ref{fig:randompredictors}b demonstrates that extending the model with $q\to\infty$ random predictors with variances $\lambda_\mathrm{opt}/q$, using the minimum-norm OLS estimator, and truncating it at $p$ dimensions is asymptotically equivalent to the ridge estimator with $\lambda_\mathrm{opt}$. As the total number of predictors $p+q$ approached $n$, MSE of the extended model increased. When $p+q$ became larger than $n$, minimum-norm shrinkage kicked in and MSE started to decrease. As $q$ grew even further, MSE approached the limiting value. In this case, $q\approx 200$ already got very close to the limiting performance.

As demonstrated in the proof, it is not necessary to truncate the minimum-norm estimator. The dashed line in Figure~\ref{fig:randompredictors}b shows the expected MSE of the full $(p+q)$-dimensional vector of regression coefficients. It converges slightly slower but to the same asymptotic value.

What if one does not know the value of $\lambda_\mathrm{opt}$ and uses random predictors with some fixed arbitrary variance to augment the model? Figure~\ref{fig:randompredictors}c shows what happens when variance is set to 1. In this case the MSE curve has a minimum at a particular $q_\mathrm{opt}$ value. This means that adding random predictors with some fixed small variance could in principle be used as an arguably bizarre but viable regularization strategy similar to ridge regression, and cross-validation could be employed to select the optimal number of random predictors.

If using random predictors as a regularization tool, one would truncate $\bbetahat_\mathrm{augm}$ at $p$ dimensions (solid line in Figures~\ref{fig:randompredictors}c). The MSE values of non-truncated $\bbetahat_\mathrm{augm}$ (dashed line) is interesting because it corresponds to the real-life $n\ll p$ situation such as in the \texttt{liver.toxicity} dataset discussed above. Our interpretation is that a small subset of high-variance PCs is actually predicting the dependent variable, while the large pool of low-variance PCs acts as an implicit regularizer.

In the simulations shown in Figure~\ref{fig:randompredictors}c, the parameter $q$ controls regularization strength and there is some optimal value $q_\mathrm{opt}$ yielding minimum expected risk. If $q<q_\mathrm{opt}$,  this regularization is too weak and some additional ridge shrinkage with $\lambda>0$ could be beneficial. But if $q>q_\mathrm{opt}$, then the regularization is too strong and no additional ridge penalty can improve the expected risk. In this situation the expected MSE as a function of $\log(\lambda)$ will be monotonically increasing on the real line, in agreement with what we saw in Figure~\ref{fig:model}d and Figure~\ref{fig:casestudy}b. Moreover, in this regime the expected MSE as a function of $\lambda$ has a minimum at a negative value $\lambda_\mathrm{opt}<0$, as we saw in Figure~\ref{fig:model}f.

We used ridge estimators on the augmented model to demonstrate this directly. Figure~\ref{fig:randompredictors}e shows the optimal ridge penalty value $\lambda_\mathrm{opt}$ for each $q$. It crosses zero around the same value of $q$ that yields the minimum risk with $\lambda=0$ (Figure~\ref{fig:randompredictors}c). For larger values of $q$, the optimal ridge penalty $\lambda_\mathrm{opt}$ is negative. This shows that negative $\lambda_\mathrm{opt}$ is due to the over-shrinkage provided by the implicit ridge regularization arising from low-variance random predictors. It is implicit over-regularization.

\subsection{Mathematical analysis for the spiked covariance model}

It would be interesting to derive some sufficient conditions on $(\boldsymbol\Sigma, \boldsymbol\beta, \sigma^2, n, p)$ that would lead to $\lambda_\mathrm{opt}\le 0$. One possible approach is to compute the derivative of $\mathbb E_{(\mathbf X, \mathbf y)} R(\bbetahat_\lambda)$ with respect to $\lambda$ at $\lambda\to 0^+$. If the derivative is positive, then $\lambda_\mathrm{opt}\le 0$.

The derivative of the risk can be computed as follows:
\begin{align}
\frac{\partial}{\partial\lambda} R (\bbetahat_\lambda) = \frac{\partial}{\partial\lambda} (\bbetahat_\lambda-\bbeta)^\top\boldsymbol\Sigma(\bbetahat_\lambda-\bbeta) = 2(\bbetahat_\lambda-\bbeta)^\top\boldsymbol\Sigma \frac{\partial\bbetahat_\lambda}{\partial\lambda}.
\end{align}
Using the standard identity  $d\mathbf A^{-1} = -\mathbf A^{-1} (d \mathbf A) \mathbf A^{-1}$, we get that
\begin{equation}
    \frac{\partial\bbetahat_\lambda}{\partial\lambda} = -(\X^\top\X + \lambda\I)^{-2}\X^\top \y.
\end{equation}
Plugging this into the derivative of the risk and setting $\lambda=0$, we obtain
\begin{equation}
    \frac{\partial}{\partial\lambda} R (\bbetahat_\lambda)\Big|_{\lambda=0} = 2\bbeta^\top\S (\X^\top\X)^{+2}\X^\top\y - 2\y^\top \X (\X^\top\X)^+ \S (\X^\top \X)^{+2}\X^\top\y,
\end{equation}
where we denote $(\X^\top\X)^{+k} = \mathbf V \mathbf S^{-2k} \mathbf V^\top$. Remembering that $\y = \X\bbeta + \boldsymbol\eta$ and taking the expectation, we get
\begin{align}
\begin{split}
    \frac{\partial}{\partial\lambda} \mathbb E_{(\X,\y)} R (\bbetahat_\lambda)\Big|_{\lambda=0} &= 2\bbeta^\top \S \mathbb E_\mathbf X (\X^\top \X)^+ \bbeta - 2\bbeta^\top \mathbb E_\mathbf X (\X^\top \X)^{+0} \S (\X^\top \X)^+ \bbeta \\
    &\quad - 2\sigma^2 \mathbb E_\mathbf X \operatorname{Tr}\Big[(\X^\top\X)^{+0}\S(\X^\top\X)^{+2}\Big],
\end{split}
\end{align}
where we used that $\mathbb E_{\boldsymbol\eta} [\mathbf a^\top \boldsymbol\eta ] = 0$ and $\mathbb E_{\boldsymbol\eta} [\boldsymbol\eta^\top \mathbf A\boldsymbol\eta] = \sigma^2 \operatorname{Tr}[\mathbf A]$ for any vector $\mathbf a$ and matrix $\mathbf A$ independent of $\boldsymbol \eta$.

We now apply this to the spiked covariance model studied above. For convenience, we write $\S = \I + c\bbeta^\top \bbeta$. Plugging this in, and denoting
\begin{equation}
    P_k = \mathbb E_\X\Big[ \bbeta^\top(\X^\top\X)^{+k}\bbeta \Big] = \mathbb E_{(\mathbf V, \mathbf S)} \Big[\bbeta^\top \mathbf V \mathbf S^{-2k} \mathbf V^\top \bbeta\Big],
\end{equation}
we obtain
\begin{equation}
    \frac{\partial}{\partial\lambda} \mathbb E_{(\X,\y)} R (\bbetahat_\lambda)\Big|_{\lambda=0} = 2c\|\bbeta\|^2P_1 - 2cP_0 P_1 - 2\sigma^2 \mathbb E_\X \operatorname{Tr}(\mathbf S^{-4}) - 2c\sigma^2 P_2. 
    \label{eq:deriv}
\end{equation}
For the spherical covariance matrix, $c=0$ and hence the derivative is always negative, in agreement with the fact that $\lambda_\mathrm{opt}>0$ for all $\bbeta$, $n$, and $p$ \citep{nakkiran2020optimal}. When $c>0$, the derivative can be positive or negative, depending on which term dominates.

We are interested in understanding the $p\gg n$ behaviour. In simulations shown above (Figures~\ref{fig:model}--\ref{fig:randompredictors}), we had  $\|\bbeta\|^2=\alpha\sigma^2/(1+\rho p) = \mathcal O(1/p)$ and $c=\rho p/\|\bbeta\|^2=\mathcal O(p^2)$. For $p\gg n$ and $0<\rho\ll 1$, all $n$ singular values of $\X$ are close to $\sqrt{p}$. This makes contribution of the third term, which is the largest when $n\approx p$ due to near-zero singular values in $\X$, asymptotically negligible because it behaves as $\mathcal O(1/p^2)$. The $\bbeta$ aligns with the leading singular vector in $\mathbf V$ and is approximately orthogonal to the others, meaning that $P_k =  \|\bbeta\|^2 \mathcal O(1/p^k)=\mathcal O(1/p^{k+1})$. Putting everything together, we see that the first, the second, and the fourth term, all behave as $\mathcal O(1/p)$.

The fourth term is roughly $\alpha/\rho$ times smaller than the first two. In our simulations $\alpha/\rho=100$, making the fourth term asymptotically negligible. The first two terms have identical asymptotic behaviour, however the first is always larger because $P_0 < \|\bbeta\|^2$. This makes the overall sum asymptotically positive, proving that $\lambda_\mathrm{opt}<0$
in the $p\to\infty$ limit.

\begin{figure}
\includegraphics[width=\textwidth]{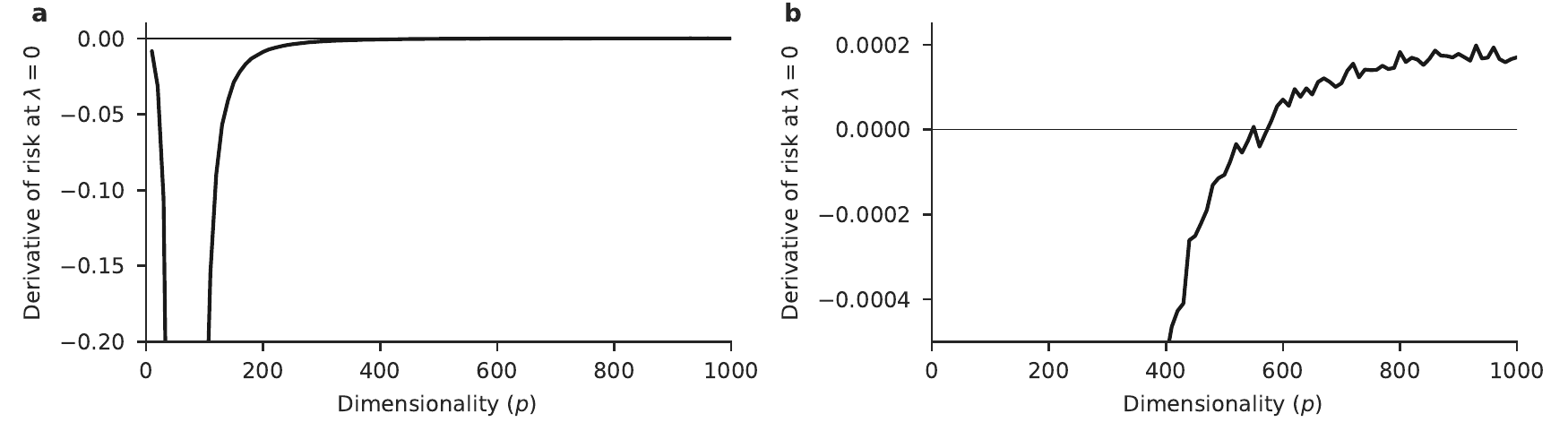}
\caption{\textbf{a.} The derivative of the expected risk as a function of ridge penalty $\lambda$ at $\lambda=0$, in the model with $p$ weakly correlated predictors. Sample size $n=64$. \textbf{b.} Zoom-in into panel (a). The derivative becomes positive for $p\gtrsim 600$, implying that $\lambda_\mathrm{opt}<0$.}
\label{fig:derivative}
\end{figure}

We numerically computed the derivative using Equation~\ref{eq:deriv} and averaging over $N_\mathrm{rep}=100$ random training set $\X$ matrices to approximate the expectation values (Figure~\ref{fig:derivative}). This confirmed that the derivative was negative and hence $\lambda_\mathrm{opt}<0$ for $p\gtrsim 600$, in agreement with Figure~\ref{fig:model}f.

\subsection{Implicit over-regularization using random Fourier features on MNIST}

For our final example, we used the setup from \citet{nakkiran2020deep, nakkiran2020optimal}, and asked whether the same phenomenon ($\lambda_\mathrm{opt}<0$) can be observed using random Fourier features on MNIST.

We normalized all pixel intensity values to lie between $-1$ and 1, and transformed the $28\times 28=784$ pixel features into 2000 random Fourier features by drawing a random matrix $\mathbf W \in \mathbb R^{784\times 1000}$ with all elements i.i.d. from $\mathcal N(0,\sigma=0.1)$, computing $\exp(-i\X\mathbf W)$ and taking its real and imaginary parts as separate features. This procedure approximates kernel regression with the Gaussian kernel, and standard deviation of the $\mathbf W$ elements corresponds to the standard deviation of the kernel \citep{rahimi2008random}. We used $n=64$ randomly selected images as a training set, and used the MNIST test set with $10\,000$ images to compute the risk. We used the digit value (from 0 to 9) as the response variable $y$, with squared error loss function. The model included the intercept which was not penalized. To estimate the expected risk, we averaged the risks over $N_\mathrm{rep}=100$ random draws of training sets.

\begin{figure}
\includegraphics[width=\textwidth]{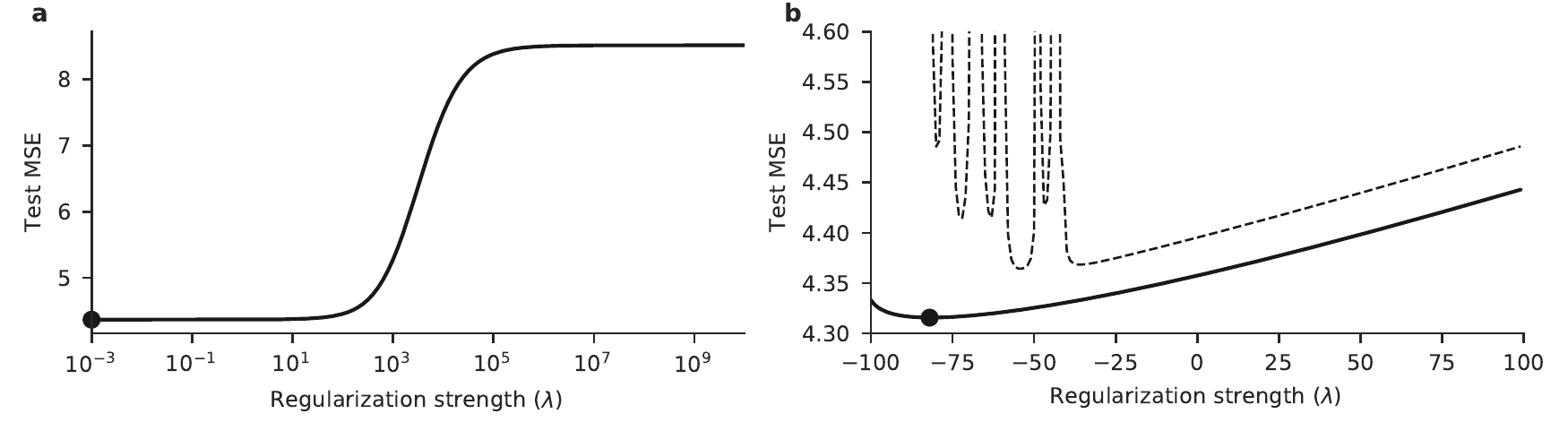}
\caption{a. Expected risk of ridge regression on MNIST data using random Fourier features as predictors and digit value as the response. Sample size $n=64$, number of Fourier features $p=2000$. When $\lambda \in \mathbb R^+$, the risk is minimized at $\lambda=0$. \textbf{b.} When $\lambda$ is allowed to take negative values, the risk is minimized at $\lambda\approx-80$, across all training sets with $s_\mathrm{min}^2>100$ (solid line; the average over 80/100 cases). Training sets with $s_\mathrm{min}^2<100$ had diverging risk around $\lambda=-s_\mathrm{min}^2$ (dashed line; the average over 20/100 cases).}
\label{fig:mnist}
\end{figure}

We found that the expected risk was minimized at $\lambda_\mathrm{opt}\approx -80$, when the expectation was computed across all 80/100 training sets that had the smallest singular value $s_\mathrm{min}^2>100$ (Figure~\ref{fig:mnist}). For any given training set, the risk diverged at $\lambda=s_\mathrm{min}$, and the smallest singular value that we observed across 100 draws was $s_\mathrm{min}^2=40$. The average risk across 20 samples with $s_\mathrm{min}^2<100$ had multiple diverging peaks for $\lambda\in[-100,0]$ (Figure~\ref{fig:mnist}b, dashed line).

The derivative of risk with respect to $\lambda$ at $\lambda=0$ that we computed in the previous section can be formally understood as the derivative at $\lambda\to0^+$. Negative derivative implies that $\lambda=0$ yields better expected risk than any positive value. However, if the generative process allows singular values of $\X$ to become arbitrarily small, then $\lambda<0$ can possibly yield diverging \textit{expected} risk. That said, for any given training set, the risk will not diverge for $\lambda\in (-s_\mathrm{min}^2, \infty)$ and the minimal \textit{conditional} risk (conditioned on the training set) can be attained at $\lambda_\mathrm{opt}<0$. Indeed, in our MNIST example, the average across all 100 training sets was monotonically decreasing until $\lambda\approx -35$ (Figure~\ref{fig:mnist}b).

\section{Discussion}

\subsubsection*{Summary and related work}

We have demonstrated that the minimum-norm OLS interpolating estimator tends to work well in the $n\ll p$ situation and that a positive ridge penalty can fail to provide any further improvement. This is because the large pool of low-variance predictors (or principal components of predictors), together with the minimum-norm requirement, can perform sufficient shrinkage on its own. This phenomenon goes against the conventional wisdom (see Introduction) but is in line with the large body of ongoing research kindled by \citet{zhang2016understanding} and mostly done in parallel to our work \citep{advani2017high, spigler2019jamming, belkin2018overfitting, belkin2018understand, belkin2019reconciling, belkin2019two, belkin2018does, nakkiran2019more, nakkiran2020deep, nakkiran2020optimal, liang2018just, hastie2019surprises, bartlett2019benign, chinot2020benign, muthukumar2019harmless, mei2019generalization, bibas2019new, derezinski2019exact}. See Introduction for more context.

We stress that the minimum-norm OLS estimator $\bbetahat_0 = \X^+ \y$ is not an exotic concept. It is given by exactly the same formula as the standard OLS estimator when the latter is written in terms of the pseudoinverse of the design matrix: $\bbetahat_\mathrm{OLS} = \X^+ \y$. When dealing with an under-determined problem, statistical software will often output the minimum-norm OLS estimator by default. 

That positive ridge penalty can fail to improve the estimator risk has been observed for kernel regression \citep{liang2018just} and for random features linear regression \citep{mei2019generalization}. Our results show that this can also happen in a simpler situation of ridge regression with Gaussian features. Our contribution is to use spiked covariance model to demonstrate and analyze this phenomenon. Moreover, we showed that the optimal ridge penalty in this situation can be negative.

In their seminal paper on ridge regression, \citet{hoerl1970ridge}  proved that there always exists some $\lambda_\mathrm{opt}>0$ that yields a lower MSE than $\lambda=0$. However, their proof was based on the assumption that $\X^\top \X$ is full rank, i.e. $n>p$. When the predictor covariance $\S$ is spherical, $\lambda_\mathrm{opt}$ is also always positive, for any $n$ and $p$ \citep{nakkiran2020deep}. Similarly, \citet{dobriban2018high} and later \citet{hastie2019surprises} proved that $\lambda_\mathrm{opt}>0$ for any $\S$ in the asymptotic $p,n\to\infty$ case while $p/n=\gamma$, based on the assumption that $\bbeta$ is randomly oriented. Here we argue that real-world $n\ll p$ problems can demonstrate qualitatively different behaviour with $\lambda_\mathrm{opt}\le 0$. This happens when $\S$ is not spherical and $\bbeta$ is pointing in its high-variance direction. This interpretation is related to the findings of \citet{bibas2019new} and \citet{bartlett2019benign}.

\subsubsection*{Augmenting the samples vs. augmenting the predictors}

It is well-known that ridge estimator can be obtained as an OLS estimator on the augmented data:
\begin{equation}\mathcal L_\lambda = \lVert \y - \X \bbeta\rVert^2 + \lambda\lVert\bbeta\rVert^2 = \bigg\lVert \begin{bmatrix}\y \\ \mathbf 0_{p\times 1}\end{bmatrix} - \begin{bmatrix}\X \\ \sqrt{\lambda}\I_{p\times p}\end{bmatrix}\bbeta\bigg\rVert^2.\end{equation}
While for this standard trick, both $\X$ and $\y$ are augmented with $p$ additional \textit{rows}, in this manuscript we considered augmenting $\X$ alone with $q$ additional \textit{columns}.

At the same time, from the above formula and from the proof of Theorem 1, we can see that if $\y$ is augmented with $q$ additional zeros and $\X$ is augmented with $q$ additional rows with all elements having zero mean and variance $\lambda/q$, then the resulting estimator will converge to $\bbetahat_\lambda$ when $q\to\infty$. This means that augmenting $\X$ with $q$ random samples and using OLS is very similar to augmenting it with $q$ random predictors and using minimum-norm OLS.

More generally, it is known that corrupting $\X$ with noise in various ways (e.g. additive noise \citep{bishop1995} or multiplicative noise \citep{srivastava2014dropout}) can be equivalent to adding the ridge penalty. Augmenting $\X$ with random predictors can also be seen as a way to corrupt $\X$ with noise.

\subsubsection*{Minimum-norm estimators in other statistical methods}

Several statistical learning methods use optimization problems similar to the minimum-norm OLS:
\begin{equation}\min \lVert\bbeta\rVert_2 \text{ subject to } \y = \X \bbeta.\end{equation}
One is the linear support vector machine classifier for linearly separable data, known to be \textit{maximum margin} classifier (here $y_i\in\{-1,1\}$) \citep{vapnik2013nature}:
\begin{equation}\min \lVert\bbeta\rVert_2 \text{ subject to } y_i (\bbeta^\top \x_i + \beta_0) \ge 1 \text{ for all }  i.\end{equation}
Another is basis pursuit \citep{chen2001atomic}:
\begin{equation}\min \lVert\bbeta\rVert_1 \text{ subject to } \y = \X \bbeta.\end{equation}

Both of them are more well-known and more widely applied in \textit{soft} versions where the constraint is relaxed to hold only approximately. In case of support vector classifiers, this corresponds to the soft-margin version applicable to non-separable datasets. In case of basis pursuit, this corresponds to basis pursuit denoising \citep{chen2001atomic}, which is equivalent to lasso \citep{tibshirani1996regression}. The Dantzig selector \citep{candes2007dantzig} also minimizes $\lVert\bbeta\rVert_1$ subject to $\y \approx \X \bbeta$, but uses $\ell_\infty$-norm approximation instead of the $\ell_2$-norm. In contrast, our manuscript considers the case where constraint $\y = \X \bbeta$ is satisfied exactly.

In the classification literature, it has been a common understanding for a long time that maximum margin linear classifier is a good choice for linearly separable problems (i.e. when $n<p$). When using hinge loss, maximum margin is equivalent to minimum norm, so from this point of view good performance of the minimum-norm OLS estimator is not unreasonable. However, when using quadratic loss as we do in this manuscript, minimum norm (for a binary $y$) is \textit{not} equivalent to maximum margin; and for a continuous $y$ the concept of margin does not apply at all. Still, the intuition remains the same: minimum norm requirement performs regularization.

\subsubsection*{Minimum-norm estimator with kernel trick}

Minimum-norm OLS estimator can be easily kernelized. Indeed, if $\x_\mathrm{test}$ is some test point, then 
\begin{equation}\hat y_\mathrm{test} = \x_\mathrm{test}^\top \bbetahat_0 = \x_\mathrm{test}^\top \X^\top (\X\X^\top)^{-1} \y = \mathbf k^\top \mathbf K^{-1} \y,\end{equation}
where $\mathbf K = \X\X^\top$ is a $n\times n$ matrix of scalar products between all training points and $\mathbf k=\X\x_\mathrm{test}$ is a vector of scalar products between all training points and the test point. The \textit{kernel trick} consists of replacing all scalar products with arbitrary kernel functions. As an example, Gaussian kernel corresponds to the effective dimensionality $p=\infty$ and so trivially $n\ll p$ for any $n$. How exactly our results extend to such $p=\infty$ situations is an interesting question beyond the scope of this paper. It has been shown that Gaussian kernel can achieve impressive accuracy on MNIST and CIFAR10 data without any explicit regularization \citep{zhang2016understanding, belkin2018understand, liang2018just} and that positive ridge regularization decreases the performance \citep{liang2018just}.

\subsubsection*{Minimum-norm estimator via gradient descent}


In the $n<p$ situation, if gradient descent is initialized at $\bbeta = 0$ then it will converge to the minimum-norm OLS solution \citep{zhang2016understanding, wilson2017marginal} (see also \citet{soudry2017implicit} and \citet{poggio2017theory} for the case of logistic loss). Indeed, each update step is proportional to $\nabla_\bbeta \mathcal L = \X^\top (\y - \X \bbeta)$ and so lies in the row space of $\X$, meaning that the final solution also has to lie in the row space of $\X$ and hence must be equal to $\bbetahat_0 = \X^+ \y = \X^\top (\X\X^\top)^{-1}\y$. If initial value of $\bbeta$ is not exactly $0$ but sufficiently close, then the gradient descent limit might be close enough to $\bbetahat_0$ to work well.

\citet{zhang2016understanding} hypothesized that this property of gradient descent can shed some light on the remarkable generalization capabilities of deep neural networks. They are routinely trained with the number of model parameters $p$ greatly exceeding $n$, meaning that such a network can be capable of perfectly fitting any training data; nevertheless, test-set performance can be very high. Moreover, increasing network size $p$ can improve test-set performance even after $p$ is large enough to ensure zero training error \citep{neyshabur2017implicit, nakkiran2020deep}, which is qualitatively similar to what we observed here. It has also been shown that in the $p\gg n$ regime, the ridge (or early stopping) regularization does not noticeably improve the generalization performance \citep{nakkiran2020optimal}.

Our work focused on \textit{why} the minimum-norm OLS estimator performs well. We confirmed its generalization ability and clarified the situations in which it can arise. Our results do not explain the case of highly nonlinear under-determined models such as deep neural networks, but perhaps can provide an inspiration for future work in that direction.

\section*{Acknowledgements}

This paper arose from the online discussion at \url{https://stats.stackexchange.com/questions/328630} in February 2018; JL and BS answered the question by DK. We thank all other participants of that discussion, in particular \texttt{@DikranMarsupial} and \texttt{@guy} for pointing out several important analogies. We thank Ryan Tibshirani for a very helpful discussion and Philipp Berens for comments and support. We thank anonymous reviewers for suggestions. DK was financially supported by the German Excellence Strategy (EXC 2064; 390727645), the Federal Ministry of Education and Research (FKZ 01GQ1601) and the National Institute of Mental Health of the National Institutes of Health under Award Number U19MH114830. The content is solely the responsibility of the authors and does not necessarily represent the official views of the National Institutes of Health.

\bibliography{main}

\end{document}